\newcommand{\be}{\begin{equation}} \newcommand{\ee}{\end{equation}}
\newcommand{\bea}{\begin{eqnarray}}\newcommand{\eea}{\end{eqnarray}}
\newcommand{\nn}{\nonumber}

\newcommand\Tq{\mbox{Tr}_q}

\documentstyle[12pt]{article}

\begin{document}
\renewcommand{\thefootnote}{\fnsymbol{footnote}}
\begin{titlepage}
December 1994\hfill{JINR-E2-94-487}
\vspace{3cm}
\begin{center}
DIFFERENTIAL CALCULUS ON THE QUANTUM SPHERE AND\\
 DEFORMED SELF-DUALITY EQUATION
\vspace{1cm}\\

B.M.Zupnik ${}\;$\footnote{E-mail: zupnik@thsun1.jinr.dubna.su}\\
Bogoliubov   Laboratory of Theoretical Physics ,\\ JINR, Dubna,
 Head Post Office,\\
P.O.Box 79, 101 000 Moscow, Russia\\
\vspace{1cm}
Talk at the International Workshop "Finite Dimensional
Integrable Systems", Dubna, 18-21 July, 1994
\vspace{2cm}\\

\end{center}
\begin{abstract}
We discuss the left-covariant 3-dimensional differential calculus on
 the quantum
sphere $SU_q (2)/U(1) $. The $SU_q (2)$-spinor harmonics are treated as
coordinates of the quantum sphere. We consider the gauge theory for
the quantum group $SU_q (2)\times U(1) $ on the deformed Euclidean space
$E_q (4)$.  A  $q$-generalization of the harmonic-gauge-field formalism
is suggested . This formalism is applied for the harmonic (twistor)
interpretation of the quantum-group self-duality equation (QGSDE).
 We consider the zero-curvature representation
and the general construction of QGSDE-solutions in terms of the analytic
prepotential.
\end{abstract}
\end{titlepage}
\renewcommand{\thefootnote}{\arabic{footnote}}
\setcounter{footnote}0
\setcounter{section}{0}
\section{ Introduction}

$\;\;\;$ The 2-dimensional sphere $S^2 $ is the simplest example of
homogeneous space and can be treated as $SU(2)/U(1) $ coset space.
 $S^2 $ plays
an important role in the twistor program of Penrose \cite{a1}
and,  particularly, in the twistor interpretation of self-duality
equation
\cite{a2}-\cite{a4}. The harmonic approach \cite{a4},\cite{a5} is a
specific
version of the twistor formalism  based on using  the spinor
harmonics as coordinates on $S^2 $.

In the present talk, we make an attempt to construct a $q$-deformed
harmonic formalism in the framework of the quantum-group concept
\cite{a6},\cite{a7}. Noncommutative geometry of  quantum spheres has
been
considered in Refs\cite{a7}-\cite{a9}. We shall use the left-invariant
3D differential calculus on the quantum group $SU_q (2)$ \cite{a10},
\cite{a11}
to study  geometry on the quantum sphere $SU_q (2)/U(1)=S_q^2 $.
Global functions on $ S_q^2 $ can be defined as the subset of
$SU_q (2)$-functions with a zero $U(1)$-charge; so we shall consider
the $SU_q (2)\times U(1) $-covariant relations for the basic
geometrical objects on $ S_q^2 $.

Quantum harmonics will be considered as matrix elements $u^i_{\pm}$
of the $SU_q (2)$-matrix $u$. An operator of external derivation $d_u$
on $SU_q (2)$ can be decomposed in terms of three invariant operators
corresponding to the different generators of a deformed Lie algebra.
We discuss the analogous decomposition of Maurer-Cartan equations on
$SU_q (2)$.

The deformed harmonic formalism can be used for analysis of the self-
duality equation on the quantum Euclidean space $E_q (4)$. The
noncommutative coordinates $x$ of $E_q (4)$ satisfy the $SU_q^L (2)
\times SU_q^R (2)$-covariant commutation relations. In this approach,
 quantum harmonics are connected with the left $SU_q (2)$-group.

We use the noncommutative algebra of differential complexes \cite{a12}-
\cite{a14} as a basis of the quantum-group gauge theory. The
 quantum-group
self-duality equation (QGSDE) on $E_q (4)$ can be formulated with the
help of a duality operation on the curvature 2-form. We present the
deformed analog of the classical BPST-instanton solution.

Quantum harmonics allow us to interpret QGSDE as a zero-curvature
equation for some harmonic decomposition of the connection form. We
discuss harmonic solutions of QGSDE by analogy with the classical
harmonic formalism \cite{a4},\cite{a5}.

\setcounter{equation}{0}
\section{ Quantum harmonics and 3D-differential calculus on the quantum
group $SU_q (2)$ }

$\;\;\;$We shall use the $R$-matrix approach \cite{a7} for definition of
 the unitary
quantum group $U_q (2) =SU_q (2)\times U(1)$ where $q$ is a real
deformation parameter. Let  $T_k^i\;\;(i,k=1,2) $ be  elements of a
quantum matrix $T$ satisfying the standard $RTT$-relations ( in the
notations of Ref\cite{a14} )
\bea
&R T T^\prime =T T^\prime  R & \\   \label{A1}
& (T)_{lm}^{ik}=T^i_l\delta^k_m ,\;\;\;(T^\prime)_{lm}^{ik}=\delta^i_l
T^k_m & \nn
\eea

The  symmetrical $R, \bar{R}\;\mbox{and}\;P^{(\pm)}$ matrices
  obey the following relations
\bea
R^2=I + \lambda R,&\;\;\;\;\bar{R}R=I,&\;\;\;\;\bar{R}=R-\lambda I\\
\label{A2}
P^{(+)} + P^{(-)}=I,&\;P^{(a)}P^{(b)} =\delta^{ab}P^{(b)},&\;R=qP^{(+)}-
q^{-1}P^{(-)} \nn
\eea
where $\lambda =q-q^{-1},\;a,b=+,-$.

It is convenient to use a covariant expression for the
$q$-generalization of an antisymmetrical symbol
\bea
&\varepsilon_{ik}(q)=\sqrt{\;q(ik)}\;\varepsilon_{ik}=-q(ik)
\varepsilon_{ki}(q) & \\ \label{A3}
 q(12)=[q(21)]^{-1}=q,&\;\;\;q(11)=q(22)=1 & \\ \nn
&\varepsilon_{ik}(q)\varepsilon^{kl}(q)=\delta_i^l\;,\;
& \label{A4}
\eea
where   $\varepsilon_{ik}$ is an ordinary
antisymmetrical symbol $(\varepsilon_{ik}=\varepsilon^{ki}) $ .

$R$-matrix elements can be written in terms of $\delta\;\mbox{and}\;
\varepsilon (q)$ symbols
\be
R^{ik}_{lm}=q\delta^i_l \delta^k_m + \varepsilon^{ki}(q)
\varepsilon_{ml}(q)    \label{A5}
\ee

 Eq(1) for the $U_q (2)$ group is equivalent to the following relations:
\bea
&\varepsilon_{ml}(q)T^l_j\;T^m_n = \varepsilon_{nj}(q)\; D(T)& \\
 \label{A6}
&\varepsilon^{ml}(q)T_l^i\;T_m^k = \varepsilon^{ki}(q)\; D(T)& \\
 \nn
\eea
where $D(T)=\mbox{Det}_q (T)$ is the quantum determinant
\be
D(T)=-\frac{q}{1+q^2}\varepsilon_{ki}(q)\varepsilon^{ml}(q)T_l^i\;T_m^k
\label{A7}
\ee

Write also the covariant relations for the inverse quantum matrix
$S(T)=T^{-1}$
\bea
&S(T^i_k)= S^i_k =\varepsilon_{kl}(q)\;T_j^l\;\varepsilon^{ji}(q)\;
 D^{-1}(T) & \\ \label{A8}
& S^i_l T^l_k = T^i_l S^l_k =\delta^i_k & \\ \nn
& T^l_i {\cal D}^m_l (q) S^k_m = {\cal D}^k_i (q)=-\varepsilon_{ji}(q)
\varepsilon^{jk}(q) &\\ \label{A9}
& S^l_i ({\cal D}^{-1})^m_l T^k_m =({\cal D}^{-1})^k_i =
-\varepsilon_{ij}(q)\varepsilon^{kj}(q) & \nn
\eea
where the notation ${\cal D}\;\mbox{and}\;{\cal D}^{-1} $ for
$SU_q (2)$-metrics is introduced.

The unitarity condition for the matrix $T$ can be formulated with the
help of involution \cite{a7}
\be
T^i_k \rightarrow \overline{T^i_k}= S^k_i \label{A10}
\ee

The condition $D(T)=1 $ corresponds to the case of $SU_q (2)$. Let us
define quantum harmonics as  matrix elements of the $SU_q (2)$-matrix
$u^i_a $. We shall distinguish the upper $SU_q (2)$ index $i=1,\;2$ and
low $U(1)$-index $a=+,\;-$. $SU_q^L (2)\times U(1)$ co-transformations
of the harmonics have the following form:
\be
u^i_{\pm} \rightarrow l^i_k u^k_{\pm}\mbox{exp}(\pm i\alpha)\label{A11}
\ee
where $\alpha$ is the $U(1)$ parameter and $l$ is the $SU_q^L (2)$-
matrix.

 Eqs(\ref{A6}) for the matrix elements $u^i_a $ are equivalent to the
basic relations
\bea
& \varepsilon_{ki}(q) u^i_{\pm}u^k_{\pm}=0 \;& \\ \nn
& \varepsilon_{ki}(q) u^i_{a}u^k_{b}=\varepsilon_{ba}(q) \;,&\\
\label{A12}
&\varepsilon^{ba}(q)u_a^i\;u_b^k = \varepsilon^{ki}(q)\;& \nn
\eea

We shall use the left-covariant 3-dimensional differential calculus
\cite{a10}, \cite{a11} for the quantum harmonics. Consider the
 $q$-traceless left-invariant 1-forms satisfying the Maurer-Cartan
equations
\bea
&\theta^a_b = \bar{u}^a_i d u^i_b & \\ \label{A13}
& \mbox{Tr}_q \theta =q\theta^+_+\;+q^{-1}\theta^-_-\;=0 &\\
 \label{A14}
& d\theta^a_b = -\theta^a_c \theta^c_b & \label{A15}
\eea
where $\bar{u}^a_i$ are components of the inverse $SU_q (2)$-harmonics.

Introduce the simple $U(1)$ notation
\be
\theta_0=\theta^+_+\;,\;\;\;\;\theta_{(+2)}=\theta^-_+\;,\;\;\;\;
\theta_{(-2)}=\theta_-^+   \label{A16}
\end{equation}

Consider the left-covariant bilinear relations between harmonics and
$\theta$-forms

\bea
& q^{\pm 2}\theta_0 u^i_{\pm}=u^i_{\pm}\theta_0 &\\ \label{A17}
& q^{\pm 1}\theta_{(p)}u^i_{\pm}= u^i_{\pm}\theta_{(p)},\;\;\;p\neq 0
  & \nn
\eea
These formulas are consistent with Eqs(\ref{A12})-(\ref{A15}). Using
the standard Leibniz rules for the operator $d$ one can obtain the
relations for the $\theta$-forms
\bea
& \theta_{(p)}^2 =0,\;\;\;\;\theta_{(+2)}\theta_{(-2)}=-q^2
\theta_{(-2)}
\theta_{(+2)} & \\ \label{A18}
& \theta_{(\pm 2)}\theta_0 =-q^{\pm 4}\theta_0 \theta_{(\pm 2)} & \nn
\eea

Consider the $SU_q (2)\times U(1)$ invariant decomposition of the
harmonic
external derivative
\bea
& d_u = \delta_0 +\delta +\bar{\delta}& \\ \label{A19}
&\delta_0 =\theta_0 D_0,\;\;\;\;\delta = \theta_{(-2)}D_{(+2)},\;\;\;\;
\bar{\delta}= \theta_{(+2)}D_{(-2)}\;& \nn
\eea
where $D_0\;\mbox{and}\;D_{(\pm 2)}$ are left-invariant differential
 operators. Note
that the $D$-operators are generators of the $q$-deformed Lie algebra
\cite{a11}
\bea
& q^2 D_{(+2)}D_{(-2)}-D_{(-2)}D_{(+2)}=D_0 & \\ \nn
& D_0 D_{(+2)}- q^4 D_{(+2)} D_0=q^2(1+q^2 )D_{(+2)}& \\ \label{A20}
& D_{(-2)}D_0-q^4 D_0 D_{(-2)}=q^2(1+q^2 )D_{(-2)}& \nn
\eea

The standard basis of the universal enveloping algebra ${\bf U}_q
[SU(2)]$
 \cite{a6} can be obtained by the nonlinear substitution \cite{a11}
\bea
& D_0=\frac{q^2}{1-q^2}(1-q^{2H} ) &\\ \label{A21}
&  D_{(\pm 2)}=q^{H/2} X^{(\pm)}& \nn
\eea

The operators $\delta_0,\;\delta\;\mbox{and}\;\bar{\delta}$ are
 nilpotent
and obey the additional condition
\be
 \left\{\delta_0,\delta\right\}+\left\{\delta_0,\bar{\delta}\right\}
+\left\{\delta,\bar{\delta}\right\}=0  \label{A22}
\ee

Define the manifest expressions for the action of these operators on
quantum harmonics
\bea
&[\delta_0, u^i_+]=u^i_+ \theta_0 ,\;\;\;\;\;
[\delta, u^i_+]=0,\;\;\;\;\;\bar{\delta} u^i_+=u^i_- \theta_{(+2)} &\\
\label{A23}
& [\delta_0, u^i_-]=- \theta_0 u^i_-\;,\;\;\;[\delta, u^i_-]=u^i_+
\theta_{(-2)},\;\;\;\;[\bar{\delta}, u^i_-]=0 & \nn
\eea

An invariant decomposition of the Maurer-Cartan equations on
$SU_q (2)/ U(1)$ has the following form:
\bea
& d_u \theta_0 =2\{\delta,\theta_0\}=2\{\bar{\delta},\theta_0\}
=-\theta_{(-2)}\theta_{(+2)} &\\ \nn
& d_u \theta_{(+2)}=2\{\delta_0,\theta_{(+2)}\}=2\{\bar{\delta},
\theta_{(+2)}\}=q^2 (1+q^2)\theta_0\theta_{(+2)}
& \\ \label{A24}
& d_u \theta_{(-2)}=2\{\delta_0,\theta_{(-2)}\}=2\{\delta,
\theta_{(-2)}\}
 =q^2 (1+q^2)\theta_{(-2)}\theta_0
& \nn
\eea

Global functions on the quantum sphere $S^2_q=SU_q (2)/ U(1)$ satisfy
the invariant condition
\be
[\delta_0,f(u)]=\theta_0 D_0 f(u)=0 \label{A25}
\ee

We shall consider also the $U(1)$-charged functions of the harmonics
$f_{(p)}(u)$
\be
 [H,f_{(p)}(u)]=pf_{(p)}(u)  \label{A26}
\ee
where $p$ is an integer number.

We shall treat harmonic functions as formal expansions on irreducible
harmonic polynomials. The $q$-symmetrized product of $r$ harmonics
$u^i_+$ and $s$ harmonics $u^i_-$ is the basis of the irreducible
 $SU_q (2)$-representation with the $U(1)$-charge $p=r-s$
\be
\Phi^{(r,s)}(u)=\Phi^{(i_1\cdots i_{r+s}) } (u)=
u^{(i_1 }_+ u^{i_2}_+\cdots u^{i_r}_+ u^{i_{r+1}}_- \cdots
u^{i_{r+s})}_- =(u_+)^r (u_-)^s  \label{A27}
\ee
where $(r,s)=I$ is the $q$-symmetrized multiindex
\be
P^{(+)}_{k,k+1} \Phi^{(r,s)}=q^{-1}R_{k,k+1} \Phi^{(r,s)}=\Phi^{(r,s)}
\label{A27b}
\ee
Here the $R$-matrix and the projectional operator $P^{(+)}$ act on the
indices
$i_k$ and  $i_{k+1}$.

The monomials $\Phi^{(r,s)}$ obey complicated commutation relations
depending on the values $r,s$ , so the polynomials $f_{(p)}(u)$
with complex numerical coefficients have not covariant commutation
properties. It is useful to extend the algebra of harmonics by adding
the set of noncommuting coefficients $C_{(r,s)}$. These  coefficients
are the components of the covariant neutral harmonic polynomials
( covariant $q$-harmonic fields )
\be
F(u)=\sum C_{(r,r)}\Phi^{(r,r)}(u)=\sum C_I \Phi^I \label{A28}
\ee

The bilinear commutation relations between $C_{I}\;\mbox{and}\;u$
follow from the requirement of harmonic commutativity :
\be
[u^i_{\pm},F(u)]=0 \label{A29}
\ee

Relations between different coefficients $C_{I}$ can be obtained , for
instance, from the additional assumption of commutativity for the
monomials in Eq(\ref{A28}). If one has a matrix harmonic field
$F^a_b (u)$
satisfying the bilinear relations, then new relations for the
corresponding
coefficients arise too.

A construction of the differential calculus on covariant harmonic
fields includes the relations for the harmonic external derivatives
(\ref{A19}) and $C_{I}$
\be
[\delta_0,C_{I}]=[\delta,C_{I}]=[\bar{\delta},C_{I}]=0 \label{A30}
\ee

\setcounter{equation}{0}
\section{ Quantum Euclidean space and quantum self-duality equation}
\indent

$\;\;\;$Quantum deformations of the Minkowski and Euclidean
4-dimensional spaces
have been considered in Refs\cite{a16}-\cite{a20}. We shall use the
coordinates
$x^i_\alpha $ of $q$-deformed Euclidean space $E_q(4)$ as generators of
 a  noncommutative algebra covariant under the coaction of the quantum
 group $G_q(4)=SU_q^L (2)\times SU_q^R (2)$
\be
x^i_\alpha \rightarrow (lxr)^i_\alpha=l^i_k  r^\beta_\alpha \otimes
x^k_\beta   \label{A31}
\ee
where $l\;\mbox{and}\;r$ are quantum matrices of the left and right
$SU_q (2)$ groups:
\bea
& R^{ik}_{lm}x^l_\alpha x^m_\beta = x^i_\gamma x^k_\rho
R^{\gamma \rho}_{\alpha \beta}& \\ \label{A32}
& R\; r\; r^\prime=r\; r^\prime\;R,\;\;\;\;\;R\;l \;l^\prime=l\;
l^\prime\;R  &\\ \label{A33}
& [r, l^\prime]=[r, x^\prime]=[l, x^\prime]=0,\;\;\;\;
\mbox{Det}_q(l)=1=\mbox{Det}_q(r) & \nn
\eea
We use two identical copies of $R$-matrices for $SU_q^L (2)\;
\mbox{and}\; SU_q^R (2)$.

The $q$-deformed central Euclidean interval $\tau$ can be constructed by
analogy with the quantum determinant
\be
\tau(x)=-\frac{q}{1+q^2}
\varepsilon^{\beta\alpha}(q)\varepsilon_{ki}(q)x^i_\alpha x^k_\beta
\label{A34}
\ee

We do not consider the quantum-group structure on $E_q (4)$ but we shall
apply the standard formula (\ref{A8}) for a definition of the
 inverse matrix $S(x)$.

It is convenient to use the following $E_q (4)$-involution:

\bea
& \overline{x^i_\alpha}=\varepsilon_{ik}(q) x^k_\beta
                         \varepsilon^{\beta\alpha}(q)
 =\tau S^\alpha_i (x) & \\ \label{A35}
& \overline{\tau}=\tau\;,\;\;\;\overline{\overline{x^i_\alpha}}=x^i_
\alpha & \nn
\eea

Let us consider the bicovariant differential calculus on the quantum
group $U_q (2)$ \cite{a21}-\cite{a24}
\bea
& T dT^\prime = R dT T^\prime R  & \\ \label{A36}
& D(T) dT = q^2 dT D(T) & \\ \label{A37}
& \omega R \omega + R\omega R \omega R =0 & \\ \label{A38}
& T \omega^\prime = R \omega R T  & \label{A39}
\eea
where $\omega^i_k(T)=dT^i_j\;S(T^j_k)$ are  the right-invariant
differential forms.

The quantum trace $\xi$ of the form $\omega$ plays an important role
in this calculus
\bea
 &\xi (T)= {\cal D}^k_i (q) \omega^i_k (T) \neq 0,\;\;\xi^2=0,\;\;
d\xi=0&  \\ \label{A40a}
& dT=\omega T=(q^2 \lambda)^{-1} [T,\xi] ,\;\;\;\; qdD(T)=\xi D(T) &\\
\label{A40b}
& d\omega=\omega^2=-(q^2 \lambda)^{-1} \left\{\xi,\omega\right\}  &
\label{A40c}
\eea

All these formulae can be used for a construction of the
$G_q (4)$-covariant differential calculus on $E_q (4)$ via the
substitution
\be
T \rightarrow x,\;\;\;\;dT \rightarrow dx,\;\;\;\;\omega(T)\rightarrow
\omega(x)=dx\;S(x)  \label{A41}
\ee

The noncommutative algebra of differential complexes \cite{a12}-
\cite{a14}
can be used for a consistent formulation of the $U_q (2)$ gauge theory
on the quantum space $E_q (4)$. Consider the $U_q (2)$ gauge matrix
$T^a_b$ defined on $E_q (4)$. Suppose that Eqs(\ref{A1},\ref{A36} -
\ref{A40c}) locally satisfy         for each "point" $x$. Coaction of
 the gauge
group $U_q (2)$ on the connection 1-form $A^a_b$ has the following
form \cite{a12}-\cite{a14}:
\bea
& A \rightarrow T(x)\;A\;S(T(x))\;+\;dT(x)\;S(T(x)) = T\;A\;S\;+\;
\omega(T)&\\ \label{A42}
& A^a_b=dx^i_\alpha A_{ib}^{\alpha a}(x) & \nn
\eea

The basic commutation relations for the form $A$ are covariant under
the gauge transformation
\be
 A \;R \;A + R \;A \;R \;A \;R\; =\;0  \label{A43}
\ee
Note  that the general relation for $A$ contains a nontrivial
 right-hand side \cite{a14}.

The restriction $\alpha=\Tq A=0 $ is inconsistent with Eq(\ref{A43}),
but we can choose the zero field-strength condition $d\alpha=\Tq dA=0 $.
This constraint for the $U(1)$-gauge field is gauge invariant.

The curvature 2-form is $q$-traceless for this model
\be
F = dA - A^2= dx^i_\alpha dx^k_\beta F^{\beta \alpha}_{ki}(x)
 \label{A44}
\ee

Basic 2-forms on $E_q (4)$ can be decomposed with the help of the
projectional operators $P^{(\pm)}$ (\ref{A2})
\bea
& dx^i_\alpha dx^k_\beta=[\;P^{(-)}dx dx^\prime P^{(+)}+
P^{(+)}dx dx^\prime P^{(-)}\;]^{ik}_{\alpha\beta}= &\\ \nn
&= \frac{q}{1+q^2}[\varepsilon^{ki}(q)d^2 x_{\alpha \beta} +
 \varepsilon_{\beta\alpha}(q)d^2 x^{ik}] & \label{A45}
\eea
By analogy with the classical case we can treat these two parts as
self-dual and anti-self-dual 2-forms under the action of a duality
operator $\ast$.

Let us consider the deformed anti-self-duality equation
\be
\ast F=-F \label{A46}
\ee

We can obtain a 5-parameter solution for the $q$-deformed anti-self-dual
$U_q (2)$-connection \cite{a24}:
\bea
&A^a_b=dx^a_\alpha\; \varepsilon_{bk}(q)\;
 \hat{x}^k_\beta \; \varepsilon^{\beta\alpha}(q) (c+\hat{\tau})^{-1}& \\
 \label{A47}
& \hat{x}^k_\beta =x^k_\beta-c^k_\beta, \;\;\;d\hat{x}=dx,\;\;\;dc=0 &\\
\nn
&R\; \hat{x}\;\hat{x}^\prime=\hat{x}\;\hat{x}^\prime\; R,\;\;\;R\;c\;
c^\prime=c\; c^\prime\; R, \;\;\;c\; x^\prime=R\;x\;c^\prime\;R & \\ \nn
& c\; dx^\prime=R\; dx\; c^\prime\; R,\;\;\;[\hat{x},\tau(\hat{x})]=0&\\
\nn
& \tau(\hat{x}) dx= q^2 dx \tau(\hat{x}) & \label{A47b}
\eea
where $c\;\mbox{and}\;c^k_\beta$ are some "parameters" and a central
function $\hat{\tau}=\tau(\hat{x})$ can be defined by substitution
$x \rightarrow \hat{x}$ in Eq(\ref{A34}).

Note that one can treat $c$ as a central periodical function which
define a solution of the first-order finite-difference equation:
$  c(\tau)=c(q^2 \tau)$. This solution is a deformed analogue of
Belavin-Polyakov-Schwarz-Tyupkin instanton. The multiparameter
$q$-generalization of the 't Hooft solution can be considered too.

\setcounter{equation}{0}
\section{ Harmonic (twistor) interpretation of$\;\;\;\;\;\;\;\;$
 quantum-group self-duality equation}

$\;\;\;$The QGSD-equation for the field strength has the following form:
\be
F^{\beta \alpha}_{ki}=[P^{(+)}FP^{(-)}]^{\alpha\beta}_{ik}=
\varepsilon_{ki}(q)F^{\beta \alpha} \label{A48}
\ee

One can obtain the integrability condition multiplying this equation
by the product of $q$-harmonics $u^i_+ u^k_+$.

Let us discuss the covariant formulation of this integrability condition
using the deformed harmonic space. It is convenient to introduce new
analytic coordinates $x_{\alpha(\pm)}$ for
$E_q(4)\otimes_q S_q^2$.
 One should use the following commutation relations

\bea
& \partial^\alpha_k x^i_\beta=\delta^\alpha_\beta \delta^i_k +
R^{ij}_{kl} R^{\alpha\rho}_{\beta\gamma}x^l_\rho \partial^\gamma_j &\\
\label{A49}
& q\partial^\alpha_i u^l_a=R_{ik}^{lm} u^k_a \partial^\alpha_m &\\
\label{A50}
& qu^i_a x^k_\beta=R^{ik}_{lm} x^l_\beta u^m_a & \\ \nn
\eea

Define the charged analytical coordinates and derivatives and
the corresponding commutation relations
\bea
& x_{\alpha a}=\varepsilon_{ab}(q)x^b_\alpha =\varepsilon_{ik}(q)
x^k_\alpha u^i_a=-q^2\varepsilon_{ki}(q)u^i_a x^k_\alpha &\\\label{A51}
& R^{cd}_{ab}x_{\alpha c}x_{\beta d}=R^{\gamma\rho}_{\alpha\beta}
x_{\gamma a}x_{\rho b}& \\ \nn
& \partial^\alpha_a =u^i_a  \partial^\alpha_i ,\;\;\;R^{cd}_{ab}
\partial^\alpha_c \partial^\beta_d =R^{\beta\alpha}_{\gamma\rho}
\partial^\rho_a \partial^\gamma_b &\\ \label{A52}
& \partial^\alpha_a x_{\beta}^b=\delta^\alpha_\beta \delta^b_a +
q^{-1}R^{\alpha\rho}_{\beta\gamma}R^{fb}_{ga}x^g_\rho \partial^\gamma_f
& \label{A53}
\eea
Note that upper and low indices $a,b\ldots$ have opposite $U(1)$-charges.

Consider the symmetrical decomposition of the external derivative $d_x$
on $E_q (4)$
\bea
&d_x=dx^i_\alpha \partial^\alpha_i=\kappa_\alpha^a \partial^\alpha_a=
d^a_a=d_1 +d_2 &\\ \label{A54}
& d_x^2 =0,\;\;\;\;\;d^a_b=\kappa_\alpha^a \partial^\alpha_b &\nn
\eea
where $\kappa_\alpha^a=\varepsilon^{ab}(q)\kappa_{\alpha b}$ are the
covariant analytic 1-forms:
\bea
& \kappa_{\alpha a}=\varepsilon_{ki}(q)dx^i_\alpha u^k_a =dx_{\alpha a}
-x_{\alpha b}\theta^b_a &\\
\label{A55}
&\{d_x,\kappa_\alpha^a \}=0,\;\;\;\;\kappa_\alpha^a \kappa_\beta^b=
-R^{ba}_{dc}\kappa_\gamma^c \kappa_\rho^d R^{\gamma \rho}_{\alpha\beta}&
\\ \label{A56}
&\partial^\alpha_a \kappa_\beta^b=R^{\alpha\rho}_{\beta\gamma}
R^{ga}_{fb}
\kappa^f_\rho \partial^\gamma_g &\nn
\eea

It is not difficult to check the following relations:
\be
d_1^2=0,\;\;\;\;\;d^2_2 +\{d_1,d_2\}=0 \label{A57}
\ee
Stress that $d^2_2\rightarrow 0$ in the limit $q\rightarrow 1$.

An analyticity condition for the functions of $x_{\alpha a}\;\mbox{and}\;
u^i_a$ has  manifest solutions $\Lambda$ depending on the analytical
coordinate $x_{\alpha (+)}$
\be
\partial^\alpha_+ \Lambda=0 \Longleftrightarrow \;d_1 \Lambda=0
\label{A58}
\ee

It should be remarked that the action of the harmonic derivatives
$\delta_0\;\mbox{and}\;\delta$ (\ref{A23}) conserves the analyticity
\be
\{\delta_0,d_1\}\Lambda=0,\;\;\;\;\;\{\delta ,d_1\}\Lambda=0 \label{A59}
\ee

Consider a decomposition of the $U_q (2)$-connection     in the central
basis (CB) (\ref{A42}) $A=a_1 + a_2$ corresponding to the decomposition
(\ref{A54}) where $a_1=\kappa_\alpha^+ A^\alpha_+ (x)$ is a connection
for the derivative $d_1$. The quantum-group self-duality equation
(\ref{A49}) is equivalent to the zero-curvature equation
\be
d_1 a_1 -a_1^2 = 0 \label{A60}
\ee

This equation has the following harmonic solution:
\be
a_1 = d_1 h\;S(h) =\omega(h, d_1 h) \label{A61}
\ee
where $h(x,u)$ is a "bridge" $U_q (2)$-matrix function. The matrix
elements of $h, d_1 h, d_x h$ and $d_u h$ satisfy the relations
analogous to  Eqs(\ref{A36}-\ref{A39}). Additional harmonic conditions
are
\be
\delta_0 h=0,\;\;\;\;d\Tq \omega(h, d h)=0 \label{A62}
\ee
where $d$ is a nilpotent operator $(d_1,\;d_x\;\mbox{or}\;d_u)$.

The bridge solution possesses a nontrivial gauge freedom
\be
h\;\to T(x)h\Lambda(x_{(+)},u)\;,\;\;\;\;\delta_0 \Lambda=d_1 \Lambda=0
\label{A63}
\ee
where $\Lambda$ is an analytical $U_q (2)$ gauge matrix.

The matrix $h$ is a transition matrix from the central basis to the
analytic basis (AB) where $d_1$ has no connection. Consider formally
the decomposition $d=d_x + d_u$ in the CB equations (\ref{A42}-\ref{A44})
although the CB-harmonic connection is equal to zero $(d_u T=0=d_u A)$.
The bridge transform is a transition to a new $u$-dependent basis
${\cal A}$ in the algebra of $U_q (2)$ differential complexes
\bea
&{\cal A}=S(h)Ah-S(h)dh=\tilde{A}_x + V &\\ \label{A64}
& \tilde{A}_x=S(h)Ah-S(h)d_x h=\kappa_\alpha^{(-)}A^\alpha_- &\\
\label{A65}
& V=v+\bar{v}=-S(h)d_u h,\;\;\;v=\theta_{(-2)}V_{(+2)},\;\;\;\;
\bar{v}=\theta_{(+2)}V_{(-2)} & \label{A66}
\eea
where $\tilde{A}_x , V, v$ and $\bar{v}$ are the AB-connection
1-forms for the operators $d_x, d_u, \delta\;\mbox{and}\;\bar{\delta}$
correspondingly.

A general solution of QGSDE can be obtained as a solution of the
basic harmonic gauge equation \cite{a4},\cite{a5}
\be
\delta h + h v =\theta_{(-2)}[D_{(+2)}h + h V_{(+2)}]=0 \label{A67}
\ee
where the connection $v$ contains the analytic prepotential $V_{(+2)}$.

We can discuss also the harmonic equations for the AB-gauge fields by
analogy with Refs\cite{a25}
\bea
&\partial^\alpha_+ V_{(+2)}=0,\;\;\;\;A^\alpha_- =
-q^{-2}\partial^\alpha_+ V_{(-2)}& \\ \label{A68}
& [D_{(+2)}+V_{(+2)}] V_{(-2)}-q^{-2}[D_{(-2)}+V_{(-2)}] V_{(+2)}=0 &
\label{A69}
\eea
where $V_{(-2)}$ is the nonanalytic gauge field for $D_{(-2)}$.

One can obtain explicit or perturbative solutions of these equations by
using the noncommutative generalizations of classical harmonic expansions
and harmonic Green functions \cite{a4},\cite{a5},\cite{a25}.

The author would like to thank V.P.Akulov, B.M.Barbashov, Ch. Devchand,
A.T. Filippov, E.A. Ivanov, J.Lukierski,
V.I.Ogievetsky, Z. Popowicz, P.N.Pyatov, A.A. Vladimirov and especially
A.P.Isaev for helpful discussions and interest in this work.

I am grateful to administration of JINR and Laboratory of Theoretical
Physics for hospitality. This work was supported in part by
International
Science Foundation (grant RUA000) and Uzbek Foundation of
 Fundamental
Researches under the contract No.40.

\end{document}